\documentclass[12pt]{amsart}
\usepackage{amscd,amssymb}
\usepackage[graph,frame,poly,arc]{xy}  
\usepackage[plainpages,backref,urlcolor=blue]{hyperref}

\theoremstyle{plain}
\newtheorem{thm}[subsection]{Theorem}

\newtheorem{prop}[subsection]{Proposition}
\newtheorem{cor}[subsection]{Corollary}

\theoremstyle{definition}
\newtheorem{rk}[subsection]{Remark}

\newtheorem{ex}[subsection]{Example}
\newtheorem{conj}[subsection]{Conjecture}

\numberwithin{equation}{section}
\newcommand{\OO}{{\mathcal O}}

\newcommand{\F}{{\mathcal F}}

\newcommand{\CC}{{\mathcal C}}

\newcommand{\al}{{\alpha}}
\newcommand{\be}{{\beta}}
\newcommand{\ep}{{\epsilon}}

\newcommand{\Z}{\mathbb{Z}}
\newcommand{\Q}{\mathbb{Q}}

\newcommand{\C}{\mathbb{C}}
\newcommand{\PP}{\mathbb{P}}

\begin{document}

\title [On free and plus-one generated curves arising from free curves]
{On free and plus-one generated curves arising from free curves by addition-deletion of a line}

\author[Alexandru Dimca]{Alexandru Dimca$^{1}$}
\address{Universit\'e C\^ ote d'Azur, CNRS, LJAD, France and Simion Stoilow Institute of Mathematics,
P.O. Box 1-764, RO-014700 Bucharest, Romania}
\email{dimca@unice.fr}

\thanks{$^1$ This work has been partially supported by the Romanian Ministry of Research and Innovation, CNCS - UEFISCDI, grant PN-III-P4-ID-PCE-2020-0029, within PNCDI III.
}



\begin{abstract} In a recent paper, after introducing the notion of plus-one generated hyperplane arrangements, Takuro Abe has shown that if we add (resp. delete) a line to (resp. from) a free line arrangement, then the resulting line arrangement is either free or plus-one generated.
In this note we prove  that the same properties hold when we replace the line arrangement  by a free curve  and add (resp. delete) a line. The proof uses a new version of a key result due originally to 
H. Schenck, H. Terao and  M. Yoshinaga, in which no quasi homogeneity assumption is needed. Two conjectures about the Tjurina number of a union of two plane curve singularities are also stated. As a geometric application, we show that, under a mild numerical condition, the projective closure of a contractible, irreducible affine plane curve is either free or plus-one generated, using a deep result due to U. Walther.
\end{abstract}
 
\maketitle


\section{Introduction} 

Let $S=\C[x,y,z]$ be the polynomial ring in three variables $x,y,z$ with complex coefficients, and let $C:f=0$ be a reduced curve of degree $d\geq 3$ in the complex projective plane $\PP^2$. 
We denote by $J_f$ the Jacobian ideal of $f$, i.e. the homogeneous ideal in $S$ spanned by the partial derivatives $f_x,f_y,f_z$ of $f$, and  by $M(f)=S/J_f$ the corresponding graded quotient ring, called the Jacobian (or Milnor) algebra of $f$.
Consider the graded $S$-module of Jacobian syzygies of $f$ or, equivalently, the module of derivations killing $f$, namely
$$D_0(f)=\{(a,b,c) \in S^3 \ : \ af_x+bf_y+cf_z=0\}.$$
According to Hilbert Syzygy Theorem, the graded Jacobian algebra $M(f)$ has a minimal resolution of the form
\begin{equation}
\label{res1}
0 \to F_3 \to F_2 \to F_1 \to F_0,
\end{equation}
where clearly $F_0=S$, $F_1=S^3(1-d)$ and the morphism $F_1 \to F_0$ is given by
$$(a,b,c) \mapsto af_x+bf_y+cf_z.$$
With this notation, the graded $S$-module of Jacobian syzygies $D_0(f)$ has the following minimal resolution
$$0 \to F_3(d-1) \to F_2(d-1).$$
We say that $C:f=0$ is an {\it $m$-syzygy curve} if the module  $F_2$ has rank $m$. Then the module $D_0(f)$ is generated by $m$ homogeneous syzygies, say $r_1,r_2,...,r_m$, of degrees $d_j=\deg r_j$ ordered such that $$1 \leq d_1\leq d_2 \leq ...\leq d_m.$$ 
We call these degrees $(d_1, \ldots, d_m)$ the {\it exponents} of the curve $C$ and $r_1,...,r_m$ a {\it minimal set of generators } for the module  $D_0(f)$. 
The smallest degree $d_1$ is sometimes denoted by $mdr(f)$ and is called the minimal degree of a Jacobian relation for $f$. 

The curve $C$ is {\it free} when $m=2$, since then  $D_0(f)$ is a free module of rank 2, see for instance \cite{KS,Sim2,ST,To}. Moreover, there are two classes of 3-syzygy curves which are intensely studied, since they are in some sense the closest to free curves.
First, we have the  {\it nearly free curves}, introduced in \cite{DStRIMS} and studied in \cite{AD, B+, Dmax,  MaVa} which are 3-syzygy curves satisfying $d_3=d_2$ and $d_1+d_2=d$. 
Then, we have the  {\it plus-one generated line arrangements} of level $d_3$, introduced by Takuro Abe in \cite{Abe} and recently studied in \cite{MP,MV},
which are 3-syzygy line arrangements satisfying $d_1+d_2=d$. In general, a  3-syzygy curve will be called a {\it plus-one generated curve}  if it satisfies $d_1+d_2=d$. { In particular, nearly free curves are special type of plus-one generated curves.}

Takuro Abe has shown that if we add (resp. we delete) a line to (resp. from) a free line arrangement in $\PP^2$, then the resulting line arrangement is either free or plus-one generated, see \cite[Theorem 1.11]{Abe}. Some of these results hold in higher dimensions as shown by Takuro Abe, see \cite[Theorem 1.4]{Abe}. 

In this note we prove the following results, showing that the same properties hold when we start with a free curve $C$ and add (resp. delete) a line. Our results can be applied to
curve arrangements with some linear components, as studied for instance in \cite{DIPS, MP,Ta}.
The first result describes which curves can be obtained by deleting a line from a free curve, containing this line as an irreducible component.
\begin{thm}
\label{thm1}
Let $C:f=0$ be a reduced curve in $\PP^2$, $L$ a line in $\PP^2$, which is not an irreducible component of $C$. We assume that the union $C'=C \cup L:f'=0$ is a free curve. Then the curve $C$ is either free or a plus-one generated curve.
\end{thm}

The second result describes which curves can be obtained by adding a line to a free curve.
\begin{thm}
\label{thm2}
Let $C:f=0$ be a reduced curve in $\PP^2$ and  $L$ a line in $\PP^2$, which is not an irreducible component of $C$. We consider the union $C'=C \cup L:f'=0$ and assume that $C$ is a free curve. Then the curve $C'$ is either free or a plus-one generated curve.
\end{thm}

In the situations considered in  Theorem \ref{thm1} and Theorem \ref{thm2}, we can describe which relations may exist between the exponents  of the  curve $C$, the exponents of the curve $C'$, the number $r=|C \cap L|$ of intersection points of $C$ and $L$ and a new invariant of the pair of curves $(C,L)$ that we introduce now. 
First we need some notation. For an isolated hypersurface singularity $(X,0)$ we set
$$\epsilon(X,0)=\mu(X,0)-\tau(X,0),$$
where $\mu(X,0)$ (resp. $\tau(X,0)$) is the Milnor (resp. Tjurina) number 
of the singularity $(X,0)$. We recall that $\epsilon(X,0) \geq 0$ and the equality holds if and only if $(X,0)$ is quasi homogeneous, see \cite{KS0}. For the curves $C_1$, $ C_2$ and $C=C_1 \cup C_2$ and a point $q \in C_1 \cap C_2$, we set
$$\epsilon(C_1,C_2)_q=\epsilon(C_1 \cup C_2,q)-\epsilon(C_1,q)$$
and then define
$$\epsilon(C_1,C_2)=\sum_{q \in C_1 \cap C_2}\epsilon(C_1,C_2)_q.$$
With this notation we have the following results.

\begin{thm}
\label{thm1B}
Let $C:f=0$ be a reduced curve in $\PP^2$, $L$ a line in $\PP^2$, which is not an irreducible component of $C$. We assume that the union $C'=C \cup L:f'=0$ is a free curve with exponents $(d_1',d_2')$. Then  the exponents $(d_1,d_2)$ (resp. $(d_1,d_2,d_3)$) of the free (resp. plus-one generated) curve $C$, the number $r=|C \cap L|$ of intersection points and the invariant $\ep=\ep(C,L)$ satisfy one of the following conditions, and all these three cases are possible.
\begin{enumerate}

\item $d_1' <d_2'$, $d_1=d_1'$ and $d_2=d_2'-1$. In this case $C$ is a free curve and $r=d_1+1-\ep$.

\item $d_1=d_1'-1$ and $d_2=d_2'$. In this case $C$ is a free curve and $r=d_2+1-\ep$.

\item $d_1=d_1'$ and $d_2=d_2'$. In this case $C$ is a plus-one generated curve and $r=d-d_3-\ep$.

\end{enumerate}

In particular, $C$ is a free curve if and only if $r\geq d_1+1-\ep$.
\end{thm}

\begin{thm}
\label{thm2B}
Let $C:f=0$ be a reduced curve in $\PP^2$ of degree $d$ and  $L$ a line in $\PP^2$, which is not an irreducible component of $C$. We consider the union $C'=C \cup L:f'=0$ and assume that $C$ is a free curve with exponents $(d_1,d_2)$. Then  the exponents $(d_1',d_2')$ (resp. $(d_1',d_2',d_3')$) of the free (resp. plus-one generated) curve $C'$, the number $r=|C \cap L|$ of intersection points and the invariant $\ep=\ep(C,L)$ satisfy one of the following conditions, and all these three cases are possible.
\begin{enumerate}

\item $d_1'=d_1$ and $d_2'=d_2+1$. In this case $C'$ is a free curve and $r=d_1+1-\ep$.

\item $d_1<d_2$, $d_1'=d_1+1$ and $d_2'=d_2$. In this case $C'$ is a free curve and $r=d_2+1-\ep$.

\item $d_1'=d_1+1$ and $d_2'=d_2+1$. In this case $C'$ is a plus-one generated curve and $r=d_3'+1-\ep$.

\end{enumerate}
In particular, $C'$ is a free curve if and only if $r\leq d_2+1-\ep$.

\end{thm}

In Example \ref{ex2B} we show that all the 6 cases in Theorems \ref{thm1B} and \ref{thm2B} can actually occur.
In Example \ref{ex2B1} we construct three (resp. two) new infinite series of free (resp. plus-one generated) conic-line arrangements which  illustrate well the cases in Theorem \ref{thm2B}.

Finally we give the following geometric application of our results above. Consider an affine plane curve $X:g(x,y)=0$ given by a reduced polynomial $g \in R=\C[x,y]$ of degree $d$. Then the projective closure
$\overline X$ of $X$ is the curve in $\PP^2$ defined by the polynomial
$$f(x,y,z)=z^dg(\frac{x}{z},\frac{y}{z}).$$
Recall that  a contractible, irreducible affine plane curve $X$ is given, in a {\it suitable global coordinate system} $(u,v)$ on $\C^2$ by the equation $u^p-v^q=0$ for some relatively prime integers $p \geq 1$ and $q \geq 1$, see \cite{GM,LZ}. In particular, $X$ has at most a unique singular point $a$, which is a cusp of type $(p,q)$, namely the singularity $(X,a)$ is given in local analytic coordinates $(u',v')$ at $a$ by the equation $u^{'p}-v^{'q}=0$. When $X$ is smooth, then $X$ is isomorphic to $\C$ and $g$ is a component of an automorphism of $\C^2$, see
\cite{Ab,Su}. We say that in this case $X$ has a cusp of type $(1,1)$.

\begin{thm}
\label{thmC}
With the above notation, assume that $X$ is irreducible and contractible and has a cusp of type $(p,q)$ such that either $p$ or $q$ is relatively prime to $d+1$, where $d= \deg X= \deg \overline X$.
Then the  projective closure $\overline X$ of the affine plane curve $X$ is either free or plus-one generated.
\end{thm}
 In view of the results in  \cite{GM,LZ}, the projective closure $\overline X$ of such a curve $X$ is a rational cuspidal curve.  Hence when $d= \deg X =\deg \overline X$ is even or $d \leq 15$ and odd, we know more precisely  that $\overline X$ is either free, or nearly free, see \cite[Theorem 3.1]{DStRIMS} and \cite[Corollary 1.3]{Mosk}.  Same result holds when $d$ is a power of a prime number, \cite[Corollary 3.2]{DStRIMS}. 
 On the other hand, for a general $d$ odd, the partial results in \cite{Mosk} do not seem to cover the result in Theorem \ref{thmC}, since they need information on the first exponent of $\overline X$. Hence 
Theorem \ref{thmC} 
 is a new support for the general conjecture, see 
\cite[Conjecture 1.1]{DStRIMS}, saying that any rational, cuspidal plane curve is either free or nearly free. 

\medskip

Concerning the invariant 
$\epsilon(C_1,C_2)_q$,
which plays a key role in our note, we propose the following.
\begin{conj}
\label{conj1}
For any two curves $C_1$ and $C_2$ without common irreducible components, and any intersection point $q \in C_1 \cap C_2$, one has
$$\epsilon(C_1,C_2)_q \geq 0.$$
\end{conj}
This conjecture clearly holds when $(C_1,q)$ is a quasi homogeneous singularity.
Note that the well known formula
$$\mu(C_1\cup C_2,q)=\mu(C_1,q)+\mu(C_2,q)+2(C_1,C_2)_q-1,$$
where $(C_1,C_2)_q$ denotes the intersection multiplicity of the curves
$C_1$ and $C_2$ at $q$, see \cite[Theorem 6.5.1]{CTC}, implies that
$$\epsilon(C_1,C_2)_q=\mu(C_1\cup C_2,q)-\mu(C_1,q)-(\tau(C_1\cup C_2,q)-\tau(C_1,q))=$$
$$=\mu(C_2,q)+2(C_1,C_2)_q-1+\tau(C_1,q) -\tau(C_1 \cup C_2,q).$$
It follows that Conjecture \ref{conj1} is equivalent to the inequality
$$\tau(C_1\cup C_2,q)\leq \tau(C_1,q)+\mu(C_2,q)+2(C_1,C_2)_q-1.$$
and it is implied by the following stronger and more symmetric  conjecture.
\begin{conj}
\label{conj2}
For any two curves $C_1$ and $C_2$ without common irreducible components, and any intersection point $q \in C_1 \cap C_2$, one has
$$\tau(C_1\cup C_2,q)\leq \tau(C_1,q)+\tau(C_2,q)+2(C_1,C_2)_q-1.$$
\end{conj}
This stronger conjecture holds when both singularities $(C_1,q)$ and $(C_2,q)$ are irreducible, as I was informed by Marcelo E. Hernandes.
\bigskip

The fact that Theorem \ref{thm1} holds when $C'$ is a conic-line arrangement was conjectured in \cite{MP}. In the same paper, the authors show that 
adding a conic to a free conic-line arrangement can yield a 4-syzygy curve, see \cite[Example 2.8]{MP}. In a subsequent paper \cite{MP2}, A. M\u acinic and P. Pokora have proven results similar to our Theorems \ref{thm1}, \ref{thm2}, \ref{thm1B} and \ref{thm2B} using a different approach, which involves however our Theorem \ref{thmB}.

I am grateful to Piotr Pokora for asking my opinion on this conjecture in \cite{MP}, which was the initial motivation for this note.

\section{Some preliminaries } 

We have the following characterizations for free and plus-one generated curves, see  \cite[Theorem 2.3]{DS}, where $D_0(f)$ is denoted by $AR(f)$.

\begin{thm}
\label{thmA}
Let $C:f=0$ be a reduced plane curve of degree $d$ and let $d_1$ and $d_2$ be the minimal degrees of a minimal system of generators for the module of Jacobian syzygies $D_0(f)$ as above.
Then the following hold.
\begin{enumerate}

\item The curve $C$ is free if and only if $d_1+d_2=d-1$.

\item The curve $C$ is plus-one generated if and only if $d_1+d_2=d$.

\item In all the other cases $d_1+d_2 >d$.

\end{enumerate}

\end{thm}

Let $I_f$ denote the saturation of the ideal $J_f$ with respect to the maximal ideal ${\bf m}=(x,y,z)$ in $S$ and consider the following  local cohomology group, usually called the Jacobian module of $f$, 
 $$N(f)=I_f/J_f=H^0_{\bf m}(M(f)).$$
We set $n(f)_k=\dim N(f)_k$ for any integer $k$ and introduce the {\it freeness defect of the curve} $C$ by the formula
$$\nu(C)=\max _j \{n(f)_j\}$$ as in \cite{AD}.
Note that $C$ is free if and only if $N(f)=0$ and hence $\nu(C)=0$, and $C$ is nearly free if and only if $\nu(C)=1$. If we set $T=3d(d-2)$, then the sequence $n(f)_k$ is symmetric with respect to the middle point $T/2$, that is one has
\begin{equation}
\label{E1}
n(f)_a=n(f)_b
\end{equation}
for any integers $a,b$ satisfying $a+b=T$, see \cite{Se, SW}. It was shown in \cite[Corollary 4.3]{DPop} that the graded $S$-module  $N(f)$ satisfies a Lefschetz type property with respect to multiplication by generic linear forms. This implies in particular the inequalities
\begin{equation}
\label{in} 
0 \leq n(f)_0 \leq n(f)_1 \leq ...\leq n(f)_{[T/2]} \geq n(f)_{[T/2]+1} \geq ...\geq n(f)_T \geq 0.
\end{equation}
Moreover, for a degree $d$ 3-syzygy curve $C$ with exponents $(d_1,d_2,d_3)$, we have the following formula for the initial degree of the graded module $N(f)$, see \cite[Theorem 3.9]{DS}.
\begin{equation}
\label{ID}
\sigma (C)= \min \{k \ : \ N(f)_k \ne 0 \}=3(d-1)-(d_1+d_2+d_3).
\end{equation}

Consider the sheafification 
$$E_C:= \widetilde{D_0(f)} $$
of the graded $S$-module $D_0(f)$, which is a rank two vector bundle on $\PP^2$, see \cite{Se} for details. Moreover, recall that
\begin{equation} \label{equa1} 
E_C=T\langle C \rangle (-1),
\end{equation}
where $T\langle C \rangle $ is the sheaf of logarithmic vector fields along $C$ as considered for instance in \cite{MaVa,Se}.
One has,  for any integer $k$, 
\begin{equation}
\label{e7}
H^0(\PP^2, E_C(k))=D_0(f)_k \text{ and }  H^1(\PP^2, E_C(k))=N(f)_{k+d-1},  
\end{equation}
where $d=\deg(f)$, for which we refer 
to \cite[Proposition 2.1]{Se}.

Consider now two reduced curves $C_1:f_1=0$ and $C_2:f_2=0$, with no common irreducible component. Let $e_1=\deg f_1$ and $e_2=\deg f_2$. We recall the following easy result, concerning the curve $C=C_1 \cup C_2: f=f_1f_2=0$, see \cite[Theorem 5.1]{DIS}.

\begin{prop}
\label{propA}
With the above notation, one has
$$mdr(f_1) \leq mdr(f) \leq mdr(f_1)+e_2.$$
\end{prop}

Next we state a new version of a key  result due to H. Schenck, H. Terao and  M. Yoshinaga, see \cite[Theorem 1.6 and Remark 1.8]{STY}. For simplicity, in this note we consider only the case $C_2$ smooth, while in the paper \cite{STY} the case $C_2$ irreducible is discussed.

\begin{thm}
\label{thmB}
With the above notation, assume that $C_2$ is a smooth curve. 
Then there is an exact  sequence of sheaves on $ \PP^2$ given by
$$ 0 \to E_{C_1}(1-e_2) \stackrel{f_2} \longrightarrow  E_C(1) \to i_{2*}\F \to 0$$
 where  $i_2: C_2  \to \PP^2$ is the inclusion and $\F=\OO_{C_2}(D)$ a line bundle  on ${C_2}$ such that
 $$\deg D=2-2g_2-r-\epsilon(C_1,C_2),$$
 where $g_{C_2}$ is the genus of the smooth curve $C_2$ and $r$ is the number of points in the reduced scheme of $C_1 \cap C_2$.
\end{thm}
In fact, in the paper \cite{STY} the authors assume that all the singularities of $C_1$, $C_2$ and $C$ are quasi homogeneous. They need this hypothesis to identify precisely the divisor $D$ with the divisor $-K_{C_2}-R$, where
$K_{C_2}$ is the canonical bundle on $C_2$ and $R$ is the divisor on $C_2$ associated to the reduced scheme of $C_1 \cap C_2$. Note that this is possible if and only if $\epsilon(C_1,C_2)=0$.
Since in our applications in this paper $C_2=L$ is a line  and a divisor on $L$ is determined by its degree,  we can work in the more general setting of Theorem \ref{thmB}.

\proof The reader will have no difficulty to check that Theorem \ref{thmB} follows easily from the arguments given in \cite{STY}. In fact, the big square diagram in \cite[Section 3]{STY}, p. 1985, shows that $\F$ is a torsion free sheaf of rank one on $C_2$. Since $C_2$ is a smooth curve, it follows that $\F$ is a line bundle, and hence has the form $\F=\OO_{C_2}(D)$  for some divisor $D$ on $C_2$. It remains to compute the degree of $D$. Then, the proof of \cite[Proposition 2.5]{STY} shows that the Hilbert polynomial of $D_0(f)$ is
\begin{equation}
\label{H1}
HP(D_0(f),t)=3 \binom{t+2}{2}-\binom{t+1+e_1+e_2}{2}+ \tau(C)=
\end{equation}
$$=3 \binom{t+2}{2}-\binom{t+1+e_1+e_2}{2}+ \mu(C)+(\tau(C)-\mu(C)),$$
where 
$$\mu(C) = \sum_{q \in C_s}\mu(C,q) \text{ and } \tau(C) = \sum_{q \in C_s}\tau(C,q),$$
with $C_s$ being the singular set of $C$. Similarly, one has
\begin{equation}
\label{H2}
HP(D_0(f_1)(-e_2),t)=3 \binom{t+2-e_2}{2}-\binom{t+1+e_1-e_2}{2}+ \tau(C_1)=
\end{equation}
$$=3 \binom{t+2-e_2}{2}-\binom{t+1+e_1-e_2}{2}+ \mu(C_1)+(\tau(C_1)-\mu(C_1),$$
where 
$$\mu(C_1) = \sum_{q \in C_{1,s}}\mu(C_1,q) \text{ and } \tau(C_1) = \sum_{q \in C_{1,s}}\tau(C_1,q),$$
with $C_{1,s}$ being the singular set of $C_1$.
It follows from the proof of  \cite[Proposition 2.5]{STY}  that the Hilbert function of the graded $S$-module given by the cokernel $\CC$ of the map
\begin{equation}
\label{H3}
0 \to D_0(f_1)(-e_2) \stackrel{f_2} \longrightarrow  D_0(f)
\end{equation}
is $$HP(\CC,t)=e_2t-r+\frac{7e_2-3e_2^2}{2}-\ep(C_1,C_2).$$
It follows that, for the shifted graded module $\CC(1)$, we have
$$HP(\CC(1),t)=e_2t-r+\frac{9e_2-3e_2^2}{2}-\ep(C_1,C_2),$$
exactly as in  \cite[Corollary  2.6]{STY}. If we set
$$\OO_{C_2}(1)=i_2^*\OO_{\PP^2}(1),$$
then one can write $\OO_{C_2}(1)=\OO_{C_2}(D')$, where the divisor $D'$ corresponds to the intersection of a line in $\PP^2$ with the curve $C_2$, and hence $\deg D'=e_2$. We recall that
$$2g_{C_2}-2=e_2^2-3e_2,$$
where $g_{C_2}$ is the genus of the smooth curve $C_2$.
Using these facts and the exact sequence in Theorem \ref{thmB}, we get, for $t$ large,
$$HP(\CC(1),t)=h^0(\OO_{C_2}(D+tD')) =\deg(D)+e_2t+\frac{3e_2-e_2^2}{2}.$$
By comparing the two formulas for $HP(\CC(1),t)$, we get the value of $\deg(D)$ as claimed in Theorem \ref{thmB}.
\endproof

\begin{rk}
\label{rkB}
We have remarked in \cite[Section 7]{DIPS}, that for the results  \cite[Theorem 1.6 and Remark 1.8]{STY} to hold,
it is enough to assume that all the singularities of $C_1$ and $C$ situated on $C_2$ are quasi homogeneous.
In fact,  it is enough to assume {\it only } that all the singularities of $C$ situated on $C_2$ are quasi homogeneous. Let $p \in C \cap C_2$ be such a singularity. Then the singularity $(C_1,p)$ is the union of all the branches of the singularity $(C,p)$, distinct from the smooth branch $(C_2,p)$. It is easy to see, using the irreducibility of $\C^*$, that a quasi homogeneous singularity $(C,p)$, say with weights
$(w_1,w_2)$, has only branches which are quasi homogeneous with exactly the same weights $(w_1,w_2)$. For more details, the reader can have a look at the proof of Proposition (7.8) in \cite{RCS}, where the relation between quasi homogeneity and $\C^*$-actions is explained.
\end{rk}

 With the above notation, by tensoring the above exact sequence with $\OO_{\PP^2}(k-1)$, for any integer $k$, we get the exact sequence
\begin{equation}
\label{e8}
0 \to E_{C_1}(k-e_2) \stackrel{f_2} \longrightarrow  E_C(k) \to i_{2*}\OO_{C_2}(D+(k-1)D') \to 0.
\end{equation}
This yields the following result by taking the the long exact sequence of cohomology associated to this exact sequence of sheaves.

\begin{cor}
\label{corB}
Consider the two reduced curves $C_1:f_1=0$ and $C_2:f_2=0$ as above and their union $C:f=0$. Assume that the curve $C_2$ is smooth.
Then there is an exact sequence  for any integer $k$ given by
$$ 0 \to D_0(f_1)_{k-e_2} \to D_0(f)_k \to H^0(C_2,\OO_{C_2}(D+(k-1)D')) \to $$
$$ \to N(f_1)_{k-e_2+e_1-1} \to N(f)_{k+e_1+e_2-1}   \to H^1(C_2,\OO_{C_2}(D+(k-1)D')), $$
with  $\deg D=2-2g_2-r-\epsilon(C_1,C_2),$
 where $g_{C_2}$ is the genus of the smooth curve $C_2$, $r$ is the number of points in the reduced scheme of $C_1 \cap C_2$ and the divisor $D'$ corresponds to the intersection of a line in $\PP^2$ with the curve $C_2$.
\end{cor}
In particular, when $C_2$ is a line $L$, then $g_2=0$ and $\deg D'=1$. Hence in this case
\begin{equation}
\label{e18}
\deg (D+(k-1)D' )=k+1-r-\epsilon(C_1,C_2).
\end{equation}

\section{The proofs of the first main results} 

\subsection{Proof of Theorem \ref{thm1}} 

Let $d=\deg C$, $d'=\deg C'=d+1$ be the degrees of these two curves.
Let $d_1' \leq d_2'$ be the exponents of the free curve $C'$, hence
$$d_1'+d_2'=d'-1=d.$$
Let $d_1 \leq d_2$ be the smallest two degrees among the degrees of a minimal set of generators for $D_0(f)$. To prove our claim it is enough to show that
$$d_1+d_2 \leq d,$$
see Theorem \ref{thmA}.
In the  conditions above, we have the following exact sequence
for  any integer $k$,
\begin{equation}
\label{EA}
 0 \to D_0(f)_{k-1}  \to D_0(f')_k \to H^0(L,\OO_L(k+1-r-\epsilon))  \to N(f)_{k+d-2}  \to 0, 
 \end{equation}
where $r$ denotes the number of points in the intersection $C \cap L$, $N(f)$ is the Jacobian module of $f$, as in Theorem \ref{thmB}, and $\epsilon=\epsilon(C,L).$

We assume that $d_1 +d_2 >d$ and show that this leads to a contradiction. There are two cases to discuss, recall Proposition \ref{propA} above.

\bigskip

\noindent  {\bf Case 1}: $d_1=d_1'-1$ and $d_2>d_2'+1$.
 
For $k<d_2'$, in the exact sequence above we have
$$\dim  D_0(f)_{k-1}  = \dim D_0(f')_k =\dim S_{k-d_1'},$$
which implies
$$\dim H^0(L,\OO_L(k+1-r-\ep)) = \dim  N(f)_{k+d-2}.$$
Hence the dimension  $\dim  N(f)_{k+d-2}$ for $k<d_2'$ is constant 
and equal to zero when $k+1-r-\ep <0$ and then increases linearly for 
$k \leq k_0=d_2'-1$. On the other hand we know that the dimensions $n_j=\dim N(f)_j$ increase  for
$j \leq j_0=3(d-2)/2$ and after that they decrease  in a symmetric way,
see the equalities \eqref{E1} and \eqref{in} above.  Note that by assumption $C$ is not free, hence $n_{j_0'} >0$, where $j_0'=[j_0]$, the integral part of $j_0$.
We have 
$$k_0+d-2=d_2'+d-3\geq j_0$$
since $d_2' \geq d/2$. When this inequality is strict, this gives already a contradiction. In the case of equality, $d$ is even, $j_0'=j_0$ and we look one step further.
For $k=d_2'$,  in the exact sequence above we have
$$\dim  D_0(f)_{k-1}  = \dim D_0(f')_k -1,$$
which implies
$n_{j_0}=n_{j_0+1}$, a contradiction with the symmetric behavior of the sequence $n_j$ in \eqref{E1}, which says that $n_{j_0-s}=n_{j_0+s}$ for any integer $s$.

\bigskip

\noindent  {\bf Case 2}: $d_1=d_1'$ and $d_2>d_2'$.

First we set $k=d_1$ and get from the exact sequence that 
$$n_{d_1+d-2}=d_1-r+1-\ep \geq 0.$$
For $d_1+1 \leq k<d_2'$, in the exact sequence above we have
$$\dim  D_0(f)_{k-1}  =\dim S_{k-1-d_1} = \binom{k+1}{2}$$
and
$$\dim  D_0(f')_{k}  =\dim S_{k-d_1} = \binom{k+2-d_1}{2}.$$
It follows that in this range for $k$ we have
$$n_{k+d-2}=\dim H^0(L,\OO_L(k+1-r-\ep))-\dim  D_0(f')_{k} +\dim  D_0(f)_{k-1}=$$
$$=k+2-r-\ep -(k+1-d_1)=d_1-r+1-\ep,$$
hence a constant value for $n_j$ when $j \leq(d_2'-1)+d-2$ . As in Case 1, we note that
$$(d_2'-1)+d-2=d_2'+d-3\geq j_0.$$
Since our curve $C$ is neither free nor nearly free, it follows that the freeness defect satisfies the inequality
$$\nu(C)=d_1-r+1-\ep \geq 2.$$
Now, for $k=d_2'$, in the exact sequence above we have
$$\dim  D_0(f)_{k-1}  =\dim S_{d_2'-d_1-1} = \binom{d_2'-d_1+1}{2}$$
and
$$\dim  D_0(f')_{k}  =\dim S_{d_2'-d_1} +\dim S_0= \binom{d_2'-d_1+2 }{2}+1.$$
It follows that  for $k=d_2'$ we have
$$n_{d_2'+d-2}=\dim H^0(L,\OO_L(d_2'+1-r-\ep))-\dim  D_0(f')_{d_2'} +\dim  D_0(f)_{d_2'-1}=$$
$$=d_2'+2-r-\ep -(d_2'-d_1+1)-1=d_1-r-\ep \geq 1.$$
This is again a contradiction with the symmetry of the sequence $n_j$ in \eqref{E1}.

\subsection{Proof of Theorem \ref{thm2}} 

Let $d=\deg C$, $d'=\deg C'=d+1$ be the degrees of these two curves.
Let $d_1 \leq d_2$ be the exponents of the free curve $C$, hence
$$d_1+d_2=d-1$$
Let $d_1' \leq d_2'$ be the smallest two degrees among the degrees of a minimal set of generators for $D_0(f')$. To prove our claim it is enough to show that
$$d_1'+d_2' \leq d',$$
see Theorem \ref{thmA}.
In the  conditions above, we have the following exact sequence
 any integer $k$,
 \begin{equation}
\label{EB}
0 \to D_0(f)_{k-1}  \to D_0(f')_k \to H^0(L,\OO_L(k+1-r-\ep))  \to  0, 
 \end{equation}
where $r$ denotes the number of points in the intersection $C \cap L$,  see Theorem \ref{thmB}, and $\epsilon=\epsilon(C,L).$
We assume that $d_1' +d_2' >d'$ and show that this leads to a contradiction. There are two cases to discuss, recall Proposition \ref{propA}.

\bigskip

\noindent  {\bf Case 1}: $d_1'=d_1+1$ and $d_2'>d_2+1$.

For $k=d_2+1$, this exact sequence gives a contradiction, since
$$\dim D_0(f)_{k-1}=\dim S_{d_2-d_1} + \dim S_0 > \dim D_0(f')_{k}=\dim S_{d_2-d_1}.$$

\bigskip

\noindent  {\bf Case 2}: $d_1'=d_1$ and $d_2'>d_2+2$.
In this case, the exact sequence yields, for $d_1 \leq k<d_2+1$, the following
$$\dim H^0(L,\OO_L(k+1-r-\ep)) =\dim D_0(f')_{k} -\dim D_0(f)_{k-1}=$$
$$=\binom{k+2-d_1}{2} -\binom{k+1-d_1} { 2}=k+1-d_1.$$
This implies the equality
$$k+1-r-\ep=k-d_1,$$
and hence $r=d_1+1-\ep$. Next, for $k=d_2+1$ we get in a similar way
from the corresponding exact sequence
$$\dim H^0(L,\OO_L(k+1-r-\ep)) =\dim D_0(f')_{k} -\dim D_0(f)_{k-1}=$$
$$=\binom{d_2+3-d_1}{ 2} -(\binom{d_2+2-d_1}{2}+1)=d_2-d_1+1.$$
But using the above value for $r$, we get
$$\dim H^0(L,\OO_L(k+1-r-\ep))=d_2+2-(d_1+1-\ep)-\ep+1=d_2-d_1+2,$$
which is a contradiction.

\subsection{Proof of Theorem \ref{thm2B}} 

The fact that these three cases are possible, even in the class of conic-line arrangements, follows from Example \ref{ex2B}. Now we prove that there are no other possibilities.
\bigskip

\noindent  {\bf Case 1} $d_1'=d_1$. Then using the exact sequence \eqref{EB} for $k=d_1-1$ we get $r>d_1-\ep$. Then, using the same sequence for $k=d_1$ we get $r=d_1+1-\ep$ and $d_2'>d_1$. By Theorem
\ref{thm2} we know that $C'$ is either free or a plus-one generated curve. In the free case, we get $d_2=d_2'+1$, which is case (1) in our Theorem. In the plus-one generated case we get $d_1'+d_2'=d'=d+1$, which implies $d_2'=d_2+2$. Now the exact sequence \eqref{EB} for 
$k=d_2$ and for $k=d_2+1$ implies that
$$\dim H^0(L,\OO_L(d_2+1-r-\ep)) =\dim H^0(L,\OO_L(d_2+2-r-\ep)) =d_2-d_1+1,$$
which is a contradiction. Therefore this situation cannot occur. 

\bigskip

\noindent  {\bf Case 2} $d_1'=d_1+1$.
If moreover $d_1=d_2$, then using the exact sequence \eqref{EB} for $k=d_1+1$ we get $d_2'=d_1'$. By Theorem
\ref{thm2} we know that $C'$ is either free or a plus-one generated curve. In the free case, using the exact sequence \eqref{EB} for all $k$, we get that $H^0(L,\OO_L(k+1-r-\ep))=0$ for all $k$, a contradiction. Therefore this situation cannot occur. It follows that $C'$ is a plus-one generated curve and the exact sequence \eqref{EB} for $k=d_3'$ shows that $r=d_3'+1-\ep$, which is case (3) in our Theorem.

Finally assume that $d_1<d_2$. If moreover $d_1'=d_2'$, we get
$d=d_1+d_2+1>2d_1+1$ and $d' \leq d_1'+d_2'=2(d_1+1)$.
This is a contradiction, since $d'=d+1$, hence this situation cannot occur.
It follows that $d_1<d_2$ and $d_1'<d_2'$. Moreover, since $C'$ is either free, or a plus-one generated curve, we have
$$d=d'-1 \leq d_1'+d_2' \leq d'=d+1.$$
It follows that $d_2=d-d_1-1 \leq d_2' \leq d-d_1=d_2+1$. If $d_2'=d_2$, this is  the case (2) in our Theorem to prove, and if $d_2'=d_2+1$ this is the case (3) in our Theorem. 

The last claim in Theorem \ref{thm2B} follows, since in the case (iii) we have $$r=d_3'+1-\ep \geq d_2'+1-\ep=d_2+2-\ep.$$

\subsection{Proof of Theorem \ref{thm1B}} 

If the curve $C$ is free, then the first two cases (i) and (ii) in 
Theorem \ref{thm1B} are just an obvious reformulation of the first two cases (i) and (ii) in Theorem \ref{thm2B}. Consider from now on the case when $C$ is a plus-one generated curve, and hence 
$$d_1+d_2=d=d_1'+d_2'.$$

\noindent  {\bf Case 1} $d_1=d_1'-1$ and $d_2=d_2'+1$.

As in the Case 1 in the proof of Theorem \ref{thm1}, we get that for 
$k <d_2'$ the dimensions $n_{k+d-2}$ is zero for $k+1-r-\ep <0$ and then they increase linearly for $r+\ep \leq k \leq k_0=d_2'-1=d_2$.
Since
$$k_0+d-2=d_2'+d-3 \geq j_0$$
as  in the proof of Theorem \ref{thm1}, the only possibility is
$d$ even and $d_2'=d/2=d_1'$. Then for $k=d_2'$ the exact sequence \eqref{EA} yields a contradiction exactly as in  Case 1 in the proof of Theorem \ref{thm1}. It follows that this case cannot occur.

\bigskip

\noindent  {\bf Case 2} $d_1=d_1'$ and $d_2=d_2'$.
We show by examples that this case may indeed occur, see Example \ref{ex2B}, case (iv).
To determine the relation between the invariants $r$, $\ep$ and the exponents of $C$, we note that in this case the formula \eqref{ID} yields
$$\sigma(C)=3(d-1)-(d+d_3)=2d-3-d_3=(d-2)+(d-1-d_3).$$
Since $k_1=d-1-d_3 <d_1$, the exact sequence \eqref{EA} for $k=k_1$ implies that
$(d-1-d_3)+1-r -\ep=0$, and hence $r=d-d_3-\ep$ as claimed.

The last claim follows from the obvious inequality
$$d-d_3 -\ep \leq d_1 -\ep.$$

\section{Examples} 

\begin{ex}
\label{ex2B}
(i) Consider $C:f=x(xy-z^2)=0$, hence $C$ is a conic union a tangent line.
Then using a computer algebra software, we see that $C$ is free with exponents $(1,1)$. If we add one more tangent line, for instance if we set
$C':f'=xy(xy-z^2)=0$ we get a free curve with exponents $(1,2)$, hence the case (1) in Theorem \ref{thm2B} really occurs. The same behavior has the curve $C'':f''=xz(xy-z^2)=0$ obtained by adding to $C$ a line passing through its singular point $(0:1:0)$.

(ii) The infinite family of conic-line arrangements constructed in \cite[Proposition 7.2]{DIPS} show that the case (2) in Theorem \ref{thm2B} really occurs.

(iii) Finally, if we add to the curve $C$ in (i) a general line, say $L:x+y=0$, we get the nearly free curve
$C':f'=x(x+y)(xy-z^2)=0$, with exponents $(2,2,2)$, and the case (3) in Theorem \ref{thm2B} really occurs. 
Similarly, if we start with the free curve $C:f=xy(xy-z^2)$ with exponents
$(1,2)$ considered above, and add a line passing through the node $(0:0:1)$, for instance $L:x+y=0$, we get a nearly free curve 
$C':f'=xy(x+y)(xy-z^2)=0$ with exponents $(2,3,3)$, showing again that 
the case (3) in Theorem \ref{thm2B} really occurs. 
To have an example with plus-one generated curves which are not nearly free, start with the free curve $C:f=(x^3+y^3)(x^3+y^3+z^3)=0$ which is free with exponents $(2,3)$. If we add a general line, say
$L: x+2y-z=0$, then the curve $C'=C \cup L$ is plus-one generated with exponents $(3,4,5)$.

(iv) In all the above examples one has $\ep=0$. Now we present two examples with $\ep>0$. If we start with the free curve $C:f=xz(xy-z^2)=0$ with exponents $(1,2)$ mentioned at point (i) of this Example, and add the line $L:x-z=0$ passing through the singular point $q=(0:1:0)$, we get a free curve $C'=C \cup L$ with exponents $(d_1',d_2')=(2,2)$, as in case (2) of Theorem \ref{thm2B}.
The singularity of $C'$ at the point $q$ is not quasi homogeneous, indeed it has the Milnor number $\mu(C',q)= 11$ and the Tjurina number $\tau(C',q)=10$. It follows that $\ep=\ep(C,L)=1$ in this case,
and $r=2=d_2+1 -\ep$ as predicted by Theorem \ref{thm2B}.

If we start with the free curve $C':f'=xz(x-z)(xy-z^2)=0$ with exponents
$(2,2)$ considered above, and delete the tangent line $x=0$,
then the new curve $C:f=z(x-z)(xy-z^2)=0$ is a plus-one generated curve, with exponents $(d_1,d_2,d_3)=(2,2,2)$, showing that the case (iii) in Theorem \ref{thm1B} really occurs. Moreover in this example
$$r=d-d_3-\ep=4-2-1=1,$$
as predicted by our result.

Assume finally that we start with the free curve  
$$C:f=xy(x-z)(y-z)(xy-z^2)=0,$$
 with exponents $(d_1,d_2)=(2,3)$. Next we add the line $L:z=0$ and get a free curve $C'=C \cup L$ with exponents $(d_1',d_2')=(3,3)$.
Hence the change from $C$ to $C'$ corresponds to case (2) in Theorem \ref{thm2B}. 
Moreover one has $r= |C \cap L|=d_2+1- \ep (C,L)=2=d_1$, and hence 
$\ep(C,L)=2$. It follows that the inequality
$r>d_1$ from \cite[Corollary 6.6 (1)]{DIS} does not hold when non quasi homogeneous singularities are involved.

\end{ex}

\begin{ex}
\label{ex2B1}
In this example we construct first a series of free conic-line arrangements. For any integer $m \geq 3$ we consider the curve
$$\CC_m: f_m=A(yz+x^2)=0$$
where $A=y \ell_2 \ldots \ell_m=0$ is the defining equation of 
 $m$ distinct lines passing through the point $q=(0:0:1)$.
The  smooth conic $Q:yz+x^2=0$ is  tangent to the line $L_1:y=0$ at the point $q$.
We show that $\CC_m$ is free with exponents $(2,m-1)$. To do this, it is enough to construct two independent syzygies $r_1$ and $r_2$ with degrees $d_1= \deg r_1=2$ and $d_2=\deg r_2=m-1$.
Note that 
$f_z=Ay$,
hence to construct a syzygy $af_x+bf_y+cf_z=0$ is equivalent to finding a linear combination $B=af_x+bf_y$ which is divisible by $Ay$. We have
\begin{equation}
\label{eqSYZ}
xf_x+yf_y=(m+1)Ayz+(m+2)Ax^2.
\end{equation}
It follows that $B=xyf_x+y^2f_y$ is divisible by $Ay$, and this yields our first syzygy $r_1$ with $d_1=2$. Next we have
$$A_yf_x-A_xf_y= A_y(A_xyz+A_xx^2+2xA)-A_x(A_y yz+Az+A_yx^2)=  2xAA_y-A_xAz.  $$
Since $A$ is divisible by $y$ but not by $y^2$, we conclude that
$A_x$ is divisible by $y$, and hence the last term $A_xAz$ is divisible by $Ay$. Note that in the term $A_y$ the only term not divisible by $y$ is
of the form $\al x^{m-1}$ with $\al \ne 0$. Hence we have to get rid of the corresponding term $2A \al x^m$ in order to get an expression $B$ divisible by $Ay$. But this can be done using equation \ref{eqSYZ}.

It follows that the total Tjurina number $\tau(\CC_m)$, that is the sum of all the Tjurina numbers of the singularities of $\CC_m$, is given by
$$\tau(\CC_m)=(m+1)^2-2(m-1)=m^2+3,$$
see \cite{Dmax}.
The curve $\CC_m$ has a complicated singularity at the point $q$, and in addition $m-1$ singularities $A_1$, coming from the intersections of the lines $L \ne L1$ with the conic. It follows that
$$\tau(\CC_m,q)=(m^2+3)-(m-1)=m^2-m+4.$$
If we add a new line $L: \ell=0$ to the curve $\CC_m$, one of the following cases may occur.
\begin{enumerate}

\item $L$ passes through the point $q$. Then the resulting curve is
of type $\CC_{m+1}$, hence free with exponents $(2,m)$. Comparing with case (1) in Theorem \ref{thmB} we see that $\ep(\CC_m,L)=1$.

\item $L$ is a tangent line to the conic $Q$, in an intersection point
$Q \cap L_j$, where $L_j: \ell_l=0$ for $2 \leq j \leq m$. Then 
$r=|\CC_m \cap L|=m$ and $\ep(\CC_m,L)=0$. For $m>3$, this corresponds to the case (2) in Theorem \ref{thmB}, hence the resulting curve $\CC_m \cup L$ is free with exponents $(3,m-1)$.

\item $L$ is the line determined by two intersection points $Q \cap L_j$ and $Q \cap L_{j'}$,  for $2 \leq j <j' \leq m$. Here $m>3$ and again
$r=|\CC_m \cap L|=m$ and $\ep(\CC_m,L)=0$. This corresponds to the case (2) in Theorem \ref{thmB}, hence the resulting curve $\CC_m \cup L$ is again free with exponents $(3,m-1)$.

\item $L$ is meeting the conic $Q$ at two points $q' \ne q''$, with
$q'=Q \cap L_j$, where $L_j: \ell_l=0$ for $2 \leq j \leq m$, and $q''$ not such an intersection point. Then $r=|\CC_m \cap L|=m+1$ and $\ep(\CC_m,L)=0$. This corresponds to the case (3) in Theorem \ref{thmB}, hence the resulting curve $\CC_m \cup L$ is a plus-one generated curve with exponents $(3,m,m)$, in other words a nearly free curve.

\item $L$ is a generic line, hence $r=|\CC_m \cap L|=m+2$ and $\ep(\CC_m,L)=0$. This corresponds to the case (3) in Theorem \ref{thmB}, hence the resulting curve $\CC_m \cup L$ is a plus-one generated curve with exponents $(3,m,m+1)$.

\end{enumerate}

\end{ex}

\begin{rk}
\label{rkS}
In Theorems \ref{thm1B} and \ref{thm2B}, the case (iii) cannot be improved in general to the claim that the new curve is nearly free.
As an example in the case of Theorem \ref{thm1B}, consider the free curve with exponents $(5,7)$ given by
$$C':f'=x(x^4+z^4)(x^8+(xz+y^2)^4)=0,$$
which is a special case of the curves considered in \cite[Remark 7.4]{DIPS}. If we delete the line $L:x=0$, we get a plus-one generated curve
$$C:f=(x^4+z^4)(x^8+(xz+y^2)^4)=0,$$
which is not nearly free, since it is plus-one generated with exponents $(5,7,9)$, as a direct computation using SINGULAR shows.
Note moreover that $|C \cap L|=2$.
To get such an  example in the case of Theorem \ref{thm2B}, it is enough to consider the case $(5)$ in Example \ref{ex2B1} above.

\end{rk}

\section{The proof of Theorem \ref{thmC}}

Let $X:g=0$ be a curve in $\C^2$, defined by a degree $d$ reduced polynomial $g$ such that in a suitable global coordinate system $(u,v)$ on $\C^2$, the curve  $X$ is given by $u^p-v^q=0$. The projective closure of $X$ is given by
$$\overline X: f(x,y,z)=z^dg(\frac{x}{z},\frac{y}{z}) =0.$$
We set $C=\overline X:f=0$ and $C'=C \cup L:f'=0$, where $L:z=0$ is the line at infinity. It follows that 
$$M=\PP^2 \setminus C' =\C^2 \setminus X.$$
The Milnor fiber of $f'$, that is the affine surface 
$$F':f'(x,y,z)-1=0$$
in $\C^3$ is the $d'$ cyclic cover of the complement $M$, where $d'=d+1=\deg f'$. More precisely, this cyclic cover corresponds to the kernel of the following map
$$\pi_1(M) \to H_1(M,\Z)=\Z \to H_1(\C^*, \Z)=\Z \to \Z/ d'\Z,$$
where the map $H_1(M,\Z)\to H_1(\C^*, \Z)$ is induced by $g$.
This composition of maps sends an elementary cycle around the irreducible curve  $X$ to $1 \in Z$. 
Let $h: \C^2 \to \C^2$ be the automorphism of $\C^2$ corresponding to the new coordinates $(u,v)$. It follows that $h(X)=X_0$, where
$$X_0: g_0=u^p-v^q=0.$$
 Hence $h$ induces an isomorphism of algebraic varieties between $M$ and $M_0=\C^2 \setminus X_0$, and hence an isomorphism of algebraic varieties between $F'$ and $F_0$, where $F_0$ is the $d'$ covering of $M_0$ associated to the morphism
$$\pi_1(M_0) \to H_1(M_0,\Z)=\Z \to \Z \to \Z/ d'\Z,$$
defined as above, with $X$ replaced by $X_0$ and $g$ replaced by $g_0$.

Consider now the surface $Z$ in $\C^3$ given by
$$Z:g_0(u,v)+w^{d'}=0$$
and the curve $D$ in $Z$ given by the intersection of $Z$ with the plane $w=0$. Note that in fact $D=X_0$. 

\medskip

\noindent  {\bf Step 1.} The Milnor fiber $F'$ is isomorphic as an algebraic variety to $Z_0=Z \setminus D$. \\
We claim in fact that 
$$ \phi: Z_0 \to M_0,  \  \  \ \phi (u,v,w)=(u,v)$$
is the $d'$ covering of $M_0$ described above. To prove this, it is enough to show that for an elementary cycle $\sigma \in \pi_1(M_0)$ about  $X_0$, the power $\sigma ^{d'}$ can be lifted to $Z_0$. This claim is obvious, since locally in the strong topology, the map $\phi$ in the neighborhood of  a smooth point of $X_0$ looks like the projection $\C \times \C^* \to \C \times \C^*$ given by $(\al,\be) \mapsto (\al,\be^{d'})$.

\medskip

\noindent  {\bf Step 2.} One has 
$$ \dim Gr^F_1H^2(F', \C)=Gr^F_1H^2(Z_0, \C)=0,$$ where $F$ denotes here the Hodge filtration on $H^2(Z_0, \C)$.\\
To prove this claim, we set $Z^*=Z \setminus \{0\}$ and $D^*=D \setminus \{0\}$, where $0$ denotes the origin of $\C^3$. Then $Z^*$ is a smooth surface, $D^*$ is a smooth divisor on $Z^*$, hence we have a Gysin  exact sequence of MHS (mixed Hodge structures), see for instance \cite[Remark C30]{STH}.
$$\ldots \to H^0(D^*,\C)(-1) \to H^2(Z^*,\C) \to H^2(Z_0,\C) \to H^1(D^*,\C) \to \ldots  .$$
Next $D^*$ is the link of the quasi homogeneous singularity $(D,0)=(X,0)$, and hence is isomorphic to $\C^*$.
It follows that the only non-zero Hodge numbers in this situation are
$$h^{1,1}(H^0(D^*,\C)(-1)) \text{ and } h^{2,2}(H^1(D^*,\C)(-1)).$$
Note also that
$$\dim Gr^F_1H^2(Z_0, \C)  =h^{1,1}(H^2(Z_0, \C))+h^{1,2}(H^2(Z_0, \C)),
$$
see for instance \cite[Appendix C]{STH} for basic properties of such Hodge numbers. Finally, $(Z,0)$ is an isolated surface singularity,
and due to the quasi homogeneity we may identify $Z^*$ with its link.
This link is a $\Q$-homology sphere since $1$ is not an eigenvalue
for the monodromy $T$ of the singularity $(Z,0)$, see for instance \cite[Proposition (3.4.7)]{STH}. Indeed, the singularity $(Z,0)$ is defined by the equation
$$u^p-v^q+w^{d+1}=0.$$
Let $\al, \be, \gamma \in \C^* \setminus \{1\}$ be such that
$$\al^p=1, \  \ \ \be ^q=1 \text{  and  } \gamma^{d+1}=1.$$
Then the eigenvalues of $T$ have the form $\lambda=\al \be \gamma$ and our claim follows from this fact, see for instance \cite[Theorem (3.4.10)]{STH}.

\medskip

\noindent  {\bf Step 3.} The curve $C'=C \cup L$ is free. \\
To prove this claim, it is enough to use the inequality
$$\dim N(f')_{d-2+k} \leq  \dim Gr^F_1H^2(F', \C)_{\gamma},$$
where $\gamma= \exp(2 \pi i/(d+1))$ and $1 \leq k \leq d+1$,
see \cite[Theorem 1.6]{Wa}. Here 
$$ Gr^F_1H^2(Z_0, \C)_{\gamma} \subset \dim Gr^F_1H^2(Z_0, \C)$$
denotes the eigenspace corresponding to the eigenvalue $\gamma$ of the monodromy operator $T'$ acting on the Milnor fiber $F'$.

Finally, we get the claim in Theorem \ref{thmC} using our Theorem \ref{thm1B}. In fact, the same proof as above gives the following similar result.

\begin{cor}
\label{corC}
Consider an affine plane curve $X:g(x,y)=0$ is given, in a  suitable global coordinate system $(u,v)$ on $\C^2$ by the equation $u^p-v^q=0$ for some relatively  integers $p \geq q \geq 1$. Assume that $d+1$ is relatively prime to both $p$ and $q$, where $d=\deg g$.
Then the  projective closure $\overline X$ of a the affine plane curve $X$, given by
$$\overline X: f(x,y,z)=z^dg(\frac{x}{z},\frac{y}{z}) =0$$
is either free or plus-one generated.
\end{cor}

\begin{rk}
\label{rkC}
In the setting of Theorem \ref{thmC}, even when the curve $C=\overline X$ is free, there is a wide range for the corresponding exponents $(d_1,d_2)$. For instance, the  curves in the family
$$C_d : f_d = z^{d-1}y + x^d + ax^2z^{d-2} + bxz^{d-1} + cz^d = 0, \  \ a \ne 0$$
for $d \geq 5$  are free with exponents $(2,d-3)$, see \cite{ST}. The corresponding affine curve
$X_d:g_d=0$ obtained by setting $z=1$ is smooth, since 
$$g_d=y+x^d+ax^2+bx+c$$
is clearly a component of the automorphism of $\C^2$ given by
$(x,y) \mapsto (x,g_d)$.

On the other hand, the curves in the family
$$C_{2k+1} : f_{2k+1} = (z^{k-1}y + x^k)^2z-x^{2k+1} = 0$$
for $k\geq 2$ are free with exponents $(k,k)$, see \cite[Theorem 4.6]{MRL}. The corresponding affine curve
$$X_{2k+1}: (y+x^k)^2-x^{2k+1}=0$$
is clearly contractible and has a cusp of type $(p,q)=(2,2k+1)$ at the origin.
\end{rk}

\begin{rk}
\label{rkC1}

Note that in the setting of Theorem \ref{thmC} one has 
$r= |C \cap L|=1,$
where $L:z=0$ is the line at infinity. This follows from the observation that 
$\PP^1$ with $s$ points deleted is contractible if and only if $s=1$.
Then Theorem \ref{thm1B} implies that
$\ep=d_1$ in case (1) and $\ep=d_2$ in case (2). As an example, for the curves $C_{2k+1}$ in Remark \ref{rkC}, we get $\ep=k$, hence the singularities at infinity of $C_{2k+1}$ and $C_{2k+1}'=C_{2k+1} \cup L$ are rather complicated. For the curve $C_{2k+1}$, this follows already from 
 \cite[Theorem 4.6]{MRL}, since this result implies that
 $$\mu(C_{2k+1},p)-\tau(C_{2k+1},p)=k^2-2k,$$
 where $p=C_{2k+1} \cap L$. Hence we get
 $$\mu(C_{2k+1}',p)-\tau(C_{2k+1}',p)=k^2-k.$$
 
\end{rk}

\end{document}